\documentclass[oneside, leqno,11pt]{amsart}
\usepackage{latexsym}
\usepackage{mathrsfs}
\usepackage{amsmath}
\usepackage{amssymb}
\usepackage[bookmarks]{hyperref}
\usepackage{pstricks,pst-plot}

\font\tenscr=rsfs10  scaled \magstep1

\font\sevenscr=rsfs7  scaled \magstep1
\font\fivescr=rsfs5  scaled \magstep1
\skewchar\tenscr='177 \skewchar\sevenscr='177 \skewchar\fivescr='177
\newfam\scrfam \textfont\scrfam=\tenscr \scriptfont\scrfam=\sevenscr
\scriptscriptfont\scrfam=\fivescr
\def\scr{\fam\scrfam}

\newcommand{\Dseq}{\{D_k\}_{k=1}^\infty}
\newcommand{\Dcoll}{{\scr D}}
\newcommand{\Kdef}{{\ol \D} \setminus \bigcup_{D\in \Dcoll} D}
\newcommand{\Dt}{\widetilde{D}}

\newcommand{\Dtcoll}{\widetilde\Dcoll}

\def\Jnseq#1{(#1)_{n=1}^\infty}
\def\Amseq#1{(#1)_{m=2}^\infty}
\def\Jsqrt#1{\sqrt{#1}}

\newcommand{\Jgood}{}

\def\Re{\mathrm{Re}}

\newtheorem{theorem}{Theorem}[section]
\newtheorem{lemma}[theorem]{Lemma}
\newtheorem{corollary}[theorem]{Corollary}
\newtheorem{proposition}[theorem]{Proposition}

\theoremstyle{definition}
\newtheorem{remark}[theorem]{Remark}
\newtheorem{definition}[theorem]{Definition}

\newcommand{\C}{\mathbb{C}}

\newcommand{\D}{\mathbb{D}}

\newcommand{\N}{\mathbb{N}}

\newcommand{\No}{\N\cup\{0\}}

\newcommand{\oD}{\overline \D}

\newcommand{\bthm}{\begin{theorem}}
\newcommand{\ethm}{\end{theorem}}
\newcommand{\blem}{\begin{lemma}}
\newcommand{\elem}{\end{lemma}}
\newcommand{\bcor}{\begin{corollary}}
\newcommand{\ecor}{\end{corollary}}
\newcommand{\bprop}{\begin{proposition}}
\newcommand{\eprop}{\end{proposition}}
\newcommand{\bdefn}{\begin{definition}}
\newcommand{\edefn}{\end{definition}}
\newcommand{\bpf}{\begin{proof}}
\newcommand{\epf}{\end{proof}}

\def \sm {\setminus}

\def\der{(d^{(k)})_{k=0}^m}

\def\dnonum {d^{(k)}}

\newcommand{\ra}{\rightarrow}

\newcommand{\ol}{\overline}

\def \Kroot{\Jsqrt{K}}
\def \Kt{\widetilde K}
\def \amk{a_m + \rho_m K_m}

\begin{document}
\title[Weakly strongly regular uniform algebras]{Weakly strongly regular uniform algebras}
\author{J. F. Feinstein}
\address{School of Mathematical Sciences, University of Nottingham, University Park, Nottingham NG7 2RD, UK }
\email{Joel.Feinstein@nottingham.ac.uk}
\author{Alexander J. Izzo}
\address{Department of Mathematics and Statistics, Bowling Green State University, Bowling Green, OH 43403}
\email{aizzo@bgsu.edu}
\thanks{This research was supported through the program \lq\lq Oberwolfach Research Fellows\rq\rq\ by the Mathematisches Forschungsinstitut Oberwolfach in 2024.   The second author was also partially supported by NSF Grant DMS-1856010.}


\begin{abstract}
Given a$\vphantom{\widehat{\widehat{\widehat{\widehat{\widehat{\widehat{\widehat X}}}}}}}$
uniform algebra $A$ on a compact Hausdorff space $X$ and a point $x$ in $X$, denote by
$M_x$ the ideal of functions in $A$ that vanish at $x$ and by $J_x$ the ideal of functions in $A$ that vanish on a neighborhood of $x$.
It is shown that for each integer $m\geq 2$, there exists a compact plane set $K$ containing the origin such that in $R(K)$ we have $\ol{J_x}\supset M_x$ for every $x\in K\sm\{0\}$ and $\ol{J_0}\supset M_0^m$ but $\ol{J_0} \not\supset M_0^{m-1}$.  This result establishes a recent conjecture of Alexander Izzo.
For the proof we introduce a construction that could be described as
taking square roots of Swiss cheeses.
\end{abstract}

\maketitle

\vskip -2.06  true in
\centerline{\footnotesize\it Dedicated to Stuart Sidney}
\vskip 2.06 true in 


%
%
%
%

\section{Introduction}

The main purpose of this paper is to prove a recent conjecture of the second author concerning closed ideals in the uniform algebra $R(K)$ for $K$ a compact plane set.  By a \emph{compact plane set} we mean a nonempty compact subset of the complex plane $\C$, and as usual, $R(K)$ denotes the uniform closure on $K$ of the holomorphic rational functions with poles off $K$.
Before presenting the conjecture, we need some background.

Let
$X$ be a compact Hausdorff space, and let $C(X)$ be the algebra of all continuous complex-valued functions on $X$ with the supremum norm
$ \|f\|_{X} = \sup\{ |f(x)| : x \in X \}$.  A \emph{uniform algebra} $A$ on $X$ is a closed subalgebra of $C(X)$ that contains the constant functions and separates
the points of $X$.
The uniform algebra $A$ on $X$ is said to be \emph{nontrivial} if $A\neq C(X)$.

Given a uniform algebra $A$ on $X$ and a point $x\in X$, we define the ideals $M_x$ and $J_x$ by
\begin{align}
M_x &= \{\, f\in A:  f(x)=0\,\}  \nonumber  \\
\intertext{and}
J_x&= \{ \, f\in A: f^{-1}(0)\ \hbox{contains a neighborhood of $x$ in $X$}\}. \nonumber
\end{align}
When it is necessary to indicate with respect to which algebra the ideals are taken,
we will denote the ideals $M_x$ and $J_x$ in the uniform algebra $A$ by $M_x(A)$ and $J_x(A)$, respectively.
The uniform algebra $A$ on $X$ is \emph{strongly regular at} $x\in X$ if $\ol{J_x}=M_x$, and $A$ is \emph{strongly regular} if $A$ is strongly regular at every point of $X$.

 It was conjectured in the 1960s that there are no nontrivial strongly regular uniform algebras.  The first counterexample was given by the first author \cite[Theorem~3.3]{F1} in 1992.  Recently the second author gave an example of the form $R(K)$ for $K$ a certain type of compact plane set called a Swiss cheese \cite[Theorem~1.2]{Izzo}.
  
Throughout the paper, $\N$ denotes the set of strictly positive integers,
$\D$ denotes the open unit disc in the complex plane, $D(a,r)$ denotes the open disc with center $a$ and radius $r$, and $\ol D(a,r)$ denotes the closure of $D(a,r)$.  Given a disc $D$, we denote the radius of $D$ by $r(D)$.

There are several different definitions in the literature for Swiss cheeses. In this paper, a \emph{Swiss cheese} is a compact subset $K$ of $\oD$ such that $K$ has empty interior and $R(K)$ is nontrivial.
The existence of such sets $K$ was first discovered by the Swiss mathematician Alice Roth \cite{Roth} and rediscovered later by S. N. Mergelyan \cite{Mer}.
Several important examples of Swiss cheeses in the literature have the form $K=\ol\D\setminus \bigcup_{k=1}^\infty D_k$ with $\Dseq$ a collection of open discs such that $\sum_{k=1}^\infty r(D_k) < 1$.
For a proof that $R(K)$ is nontrivial for such $K$ see
\cite[Lemma~1]{Mc} or \cite[Lemma~24.1]{Stout} or for an indirect proof via \lq classicalisation' of Swiss cheeses \cite[Theorem 2.1]{FH}. In this paper, however, nontriviality will follow directly from other properties of $R(K)$.

The construction of the strongly regular $R(K)$ in \cite{Izzo} involved first producing a Swiss cheese $K_1$ such that in $R(K_1)$ the inclusion $\ol{J_x}\supset {M_x^2}$ holds for every $x\in K_1$, and then combining this with John Wermer's construction \cite{W1} of a Swiss cheese $K_2$ such that in $R(K_2)$ the equality $\ol{M_x^2}=M_x$ holds for every $x\in K_2$.
The question thus arose as to whether for an $R(K)$ (with $K$ a compact plane set) either of the conditions $\ol{M_x^2}=M_x$ for every $x\in K$ and $\ol{J_x}\supset {M_x^2}$ for every $x\in K$ implies the other.  
The second author gave an example showing that the first condition does not imply the second \cite[Theorem~1.4]{Izzo2}, and he conjectured that the second condition does not imply the first.  More generally, he conjectured that, for each integer $m\geq 2$, the condition $\ol{J_x}\supset {M_x^m}$  for every $x\in K$ does not imply that $\ol{J_x}\supset {M_x^{m-1}}$ (and hence does not imply that $\ol{M_x^m}=\ol{M_x^{m-1}}$) for every $x\in K$ 
\cite[Conjecture~1.14]{Izzo}.  The main purpose of the present paper is to prove
the following theorem which is slightly stronger than that conjecture.

\bthm\label{main-theorem}
For each integer $m\geq 2$, there exists a Swiss cheese $K$ containing the origin such that $R(K)$ is normal and in $R(K)$ we have $\ol{J_x}\supset M_x$ for every $x\in K\sm\{0\}$ and $\ol{J_0}\supset M_0^m$ but $\ol{J_0}\not\supset M_0^{m-1}$.
\ethm

For $A$ a uniform algebra on a compact space $X$ and $x\in X$,
we will say that $A$ is \emph{weakly strongly regular at} $x$ if there exists $m \in \N$ such that $\ol{J_x}\supset M_x^m$, and that the uniform algebra $A$ is \emph{weakly strongly regular} if $A$ is weakly strongly regular at every point of $X$.
This terminology is analogous to the definition of weak spectral synthesis for general commutative Banach algebras used in
\cite{KU,Wa}.
Note that with this terminology, for each of the compact plane sets $K$ constructed in Theorem \ref{main-theorem}, $R(K)$ is weakly strongly regular but is not strongly regular.

 A uniform algebra $A$ on a compact space $X$ is said to be
\begin{enumerate}
\item[(a)] \emph{natural} if the maximal ideal space of $A$ is $X$ (under the usual identification of a point of $X$ with the corresponding complex homomorphism),
\item[(b)] \emph{regular on $X$} if for each closed subset $K_0$ of $X$ and each point $x$ of $X\setminus K_0$, there exists a function $f$ in $A$ such that $f(x)=1$ and $f=0$ on $K_0$,
\item[(c)] \emph{normal on $X$} if for each pair of disjoint closed subsets $K_0$ and $K_1$ of $X$, there exists a function $f$ in $A$ such that $f=1$ on $K_1$ and $f=0$ on $K_0$.
\end{enumerate}
The uniform algebra $A$ on $X$ is \emph{regular} or \emph{normal} if $A$ is natural and is regular on $X$ or normal on $X$, respectively.

It is easily shown that every weakly strongly regular uniform algebra is natural by repeating an argument given by Raymond Mortini~\cite[Proposition~2.4]{Mortini} to reprove Wilken's result~\cite[Lemma]{Wilken} that every strongly regular uniform algebra is natural.  
This was done in \cite[Theorem~1.10]{Izzo} under the stronger hypothesis that there is a fixed positive integer $m$ such that $\ol{J_x}\supset M_x^m$ for every $x$, but relaxing the hypothesis to allow $m$ to depend on $x$ has no effect on the argument.  
Since it is easily shown that a weakly strongly regular uniform algebra $A$ on $X$ is regular on $X$, and since every uniform algebra that is regular on its maximal ideal space is normal \cite[Theorem~27.2]{Stout}, it follows that every weakly strongly regular uniform algebra is normal.  In particular, the assertion in the conclusion of Theorem~\ref{main-theorem} that $R(K)$ is normal is redundant.  (Note also that it is standard that, for $K$ a compact plane set,  $R(K)$ is always natural.)

The question arises as to whether to each weakly strongly regular uniform algebra $A$ there corresponds a fixed positive integer $m$ such that $\ol{J_x}\supset M_x^m$ for every point $x$ in the maximal ideal space of $A$.  As a corollary of Theorem~\ref{main-theorem}, we will show that this is not the case.

\bcor\label{m-to-infinity}
There exists a compact plane set $K$ such that $R(K)$ is weakly strongly regular and such that for every integer $m\geq 2$ there exists a point $x\in K$ such that $\ol{J_x}\supset M_x^m$ but $\ol{J_x}\not\supset M_x^{m-1}$.

\ecor

The first author \cite[Theorem~5.1]{F0} constructed a Swiss cheese $K$ that contains the origin and is such that $R(K)$ has the property that $\ol{J_0}\supset M_0^2$ but $\ol{J_0} \not\supset M_0$.  The condition that $\ol{J_0}\supset M_0^m$ but $\ol{J_0} \not\supset M_0^{m-1}$ in Theorem~\ref{main-theorem} will be obtained by starting with that Swiss cheese and successively applying a procedure that, very roughly, amounts to taking a square root of a Swiss cheese.  (See Theorems~\ref{square-root-of-set} and~\ref{square-root-of-cheese} below.)

We omit the elementary proof of the following well-known, easy lemma.

\blem\label{straightforward}
Given $m\in\N$ and $M$ an ideal in a commutative Banach algebra $A$, if $\ol{M^m}=\ol{M^{m+1}}$, then $\ol{M^m}= \ol{M^q}$ for all $q\geq m$.
\elem

\medskip
The standard fact that every ideal in a  normal uniform algebra is local \cite[Proposition~4.1.20(iv)]{D} yields at once that in every normal uniform algebra we have that $\ol{M_x^j} \supset \ol{J_x}$ for every positive integer $j$.  Thus by Lemma~\ref{straightforward}, we have in Theorem~\ref{main-theorem} that
$M_0\supsetneq \ol{M_0^2} \supsetneq \cdots \supsetneq \ol{M_0^m} =\ol{M_0^{m+1}} = \cdots$, that is, the sequence of closed powers of $M_0$ strictly decreases to $\ol{M_0^m}$ ($=\ol{J_0}$) and then stabilizes.  This phenomenon that the sequence of closed powers of a maximal ideal in a uniform algebra can stabilize at any point in the sequence was noted long ago by Stuart Sidney who studied the behavior of the sequence of closed powers of a maximal ideal much more generally \cite{Sidney}.  However, the uniform algebras constructed by Sidney were not of the form $R(K)$ for a compact plane set $K$.  Furthermore, he did not deal with the ideal $\ol{ J_x}$, whereas for us, the relation between $\ol{J_x}$ and the closed powers of $M_x$ is the primary concern.  Actually, that the sequence of closed powers of a maximal ideal in a uniform algebra can stabilize at any point in the sequence is quite easy to prove; it follows immediately from the easy Lemma~\ref{road-runner} below and a lemma about the relationship between the closed powers of a maximal ideal and higher-order bounded point derivations (Lemma~\ref{derivations-vs-ideals}).  (Higher-order point derivations on a uniform algebra are defined in the next section.)

A useful criterion for the existence of a nonzero bounded point derivation on $R(K)$ was found by Andrew Browder and presented in~\cite{W1}.  We will use a generalization of a modification of that condition to investigate higher-order point derivations.

Given a disc $D\subset\C$ and a point $a\in \C$, we denote the distance from $\partial D$ to $a$ by $s_a(D)$.  Explicitly, $s_a(D)= \inf \{ |z-a| : z\in \partial D\}$.

Let $\Dcoll$ be a \Jgood collection of open discs, let $m \in \No$, and let
$a \in \Kdef$.
We define the \emph{$m$th order Browder sum for $\Dcoll$ at $a$} to be
\[
\sum_{D\in \Dcoll \cup \{\D\}} \frac{r(D)}{s_a(D)^{m+1}}\,
\]
and we say that $\Dcoll$ \emph{satisfies the Browder condition of order $m$ at $a$} if
this Browder sum is finite.  We say that a compact plane set $K$ \emph{satisfies the Browder condition of order $m$ at $a$} if $K=\Kdef$ for some \Jgood collection of open discs $\Dcoll$ that satisfies the Browder condition of order $m$ at $a$.
For us the object of interest will always be the compact plane set $K=\Kdef$ rather than the collection $\Dcoll$.  \emph{Since the set $K$ is unaffected by any discs in $\Dcoll$ that do not intersect $\D$, we can always assume without loss of generality, and therefore we will always tacitly assume, that every disc in $\Dcoll$ intersects $\D$.}
Then $s_a(D)\leq 2$ for every $a\in \D$ and $D\in\Dcoll$, so
the Browder condition of order $m$ implies the Browder condition of order $k$ for all $k\leq m$.

In the next section we will see why the following lemma is essentially immediate from \cite[Lemma 4.4]{FY}.

\blem\label{Browder}
Let $\Dcoll$ be a \Jgood collection of open discs, set
$K= \Kdef$,
and suppose that $a \in K$.
If $\Dcoll$ satisfies the Browder condition of order $m$ at $a$ for some $m\in\N$, then there is a nondegenerate bounded point derivation of order $m$ on $R(K)$ at $a$. If $\Dcoll$ satisfies the Browder condition of order $m$ at $a$ for every $m\in \N$, then there is a nondegenerate bounded point derivation of infinite order on $R(K)$ at $a$.
\elem

A Swiss cheese $K$ such that $R(K)$ is normal but has a nondegenerate bounded point derivation of infinite order at the origin was constructed by Anthony O'Farrell~\cite{OF2}.
As a step toward proving Theorem~\ref{main-theorem}, we will strengthen O'Farrell's result as follows.

\bthm\label{infinite-order}
There exists a compact plane set $K$ containing the origin such that $R(K)$ is normal, $R(K)$ is strongly regular at every point other than the origin, and $R(K)$ has a nondegenerate bounded point derivation of infinite order at the origin.  Moreover, $K$ can be taken to be a Swiss cheese $\Kdef$ where $\Dcoll$ is a \Jgood collection of open discs which satisfies the Browder condition of order $m$ at the origin for every $m\in \N$.
\ethm

\medskip
In \cite{F2} the first author proved that regularity of $R(K)$ does not pass to finite unions.
We will obtain the following result showing that the same holds for strong regularity
as an easy consequence of Theorem~\ref{infinite-order}.

\bcor\label{unions}
There exist two compact plane sets $K_1$ and $K_2$ such that each of $R(K_1)$ and $R(K_2)$ is strongly regular, but $R(K_1\cup K_2)$ is not strongly regular.  Furthermore, $K_1$ and $K_2$ can be chosen in such a way that, in addition, $R(K_1\cup K_2)$ is normal and is strongly regular at every point of $K_1\cup K_2$ with one exception where $R(K_1\cup K_2)$ has a nondegenerate bounded point derivation of infinite order.
\ecor

For a subset $E$ of $\C$ we define the \emph{square root $\Jsqrt E$ of $E$} to be the set
$\{z \in \C: z^2 \in E\,\}$.
For the proof of Theorem~\ref{main-theorem} we will need, for each integer $m$, a Swiss cheese $K$ that satisfies the Browder condition of order $m$ at the origin and has the property that $\ol{J_0(R(K))} \supset {M_0(R(K))^q}$ for some $q$.  As alluded to earlier, this we will obtain by taking square roots of Swiss cheeses.  Specifically, we will prove the following theorems.

\bthm\label{square-root-of-set}
Let $K$ be a compact plane set that contains the origin.
\begin{enumerate}
\item[(a)] If $z^m \in \ol{J_0(R(K))}$, then $z^{2m}\in \ol{J_0(R(\Kroot))}$.
\item[(b)] If $\ol{J_0(R(K))} \supset {M_0(R(K))}^m$, then
$\ol{J_0(R(\Kroot))} \supset {M_0(R(\Kroot))}^{2m}$.
\item[(c)] If $R(K)$ has a nondegenerate bounded point derivation of order $m$ at the origin, then $R(\Kroot))$ has a nondegenerate bounded point derivation of order $2m$ at the origin.
\end{enumerate}
\ethm

\bthm\label{square-root-of-cheese}
Let $\Dcoll$ be a \Jgood collection of open discs, set
$K= \Kdef$,
and suppose that $0 \in K$. Let $m\in\N$, and suppose that $\Dcoll$ satisfies the Browder condition of order $m$ at $0$. Then there exists a \Jgood collection $\Dtcoll$ of open discs such that the following properties hold:

\begin{enumerate}
\item[(a)] $0 \in \ol\D \setminus \bigcup_{\Dt\in \Dtcoll} \Dt \subseteq \Kroot$;
\item[(b)] $\Dtcoll$ satisfies the Browder condition of order $2m$ at $0$;
\item[(c)]  the $2m$th order Browder sum for $\Dtcoll$ at $0$ is strictly less than the $m$th order Browder sum for $\Dcoll$ at $0$ (provided $\Dcoll$ is nonempty).
\end{enumerate}
\ethm

In the next section we define higher-order point derivations and present results we will need about them.  In Section~\ref{Square roots of Swiss cheeses} we prove Theorems~\ref{square-root-of-set} and~\ref{square-root-of-cheese} concerning square roots of Swiss cheeses.  Section~\ref{Some preliminary results} gives further results we will need  for the proofs of  Theorems~\ref{main-theorem} and~\ref{infinite-order} and Corollaries~\ref{m-to-infinity} and~\ref{unions}, which are given in Section~\ref{The Proofs of the main theorems}.

It is a pleasure to dedicate this paper to Stuart Sidney with whom the authors have always enjoyed talking.  As noted above, Sidney's work on the sequence of closed powers of a maximal ideal \cite{Sidney} is related to the results of our paper.  In addition, we each found Sidney's approach to nonlocal uniform algebras helpful with other projects.

%
%
%
%

\section{Higher-order point derivations and ideals}\label{preliminaries}

\medskip
For $\phi$ a point of the maximal ideal space of a uniform algebra $A$, a \emph{point derivation} on $A$ at $\phi$ is a linear functional $d$ on $A$ satisfying the identity
\begin{equation}\label{derivation-eq}
d(fg)=d(f)\phi(g) + \phi(f)d(g)
\end{equation}
for all $f$ and $g$ in $A$.
A point derivation is said to be \emph{bounded} if it is bounded (continuous) as a linear functional.
Now let $m$ be a positive integer or $\infty$.  A \emph{point derivation of order $m$ at $\phi$} is a sequence $d=\der$ of linear functionals on $A$ such that $d^{(0)}=\phi$ and for all $f$ and $g$ in $A$
\begin{equation}\label{k}
\dnonum(fg)=\sum_{j=0}^k (d^{(j)}f) (d^{(k-j)}g)\quad \hbox{for $1\leq k\leq m\ ({\rm or}\ 1\leq k<\infty)$}.
\end{equation}
The point derivation $\der$ is \emph{bounded} if each $d^{(k)}$ is bounded; $\der$ is \emph{nondegenerate} if $d^{(1)}\neq 0$.  When $\der$ is nondegenerate the functionals $d^{(0)}, d^{(1)}, d^{(2)}, \ldots$ are linearly independent.

Let $K$ be a compact plane set.
We denote by $R_0(K)$ the dense subalgebra of $R(K)$ consisting of the rational functions with poles off $K$. When discussing continuity, we give $R_0(K)$ the uniform norm on $K$.
For each $a\in K$ and each $m\in \No$, we define the linear functional $\delta_{a,m}$ on $R_0(K)$ to be the map $f \mapsto f^{(m)}(a)/m!$. Thus, for $f $ a rational function, 
$\delta_{a,m}(f)$ is the coefficient of $(z-a)^m$ in the Taylor series for $f$ at $a$.
It is standard and easily verified that there is a nonzero bounded point derivation on $R(K)$ at $a$ if and only the linear functional $\delta_{a,1}$ is bounded on $R_0(K)$.

The following is essentially \cite[Lemma 4.4]{FY} stated in terms of Browder sums and the functionals $\delta_{a,m}$. For an explicit proof of the estimate for the first derivative, see the proof of \cite[Lemma 2.11]{F0}.

\blem\label{Browder2}
Let $\Dcoll$ be a \Jgood collection of open discs, set
$K= \Kdef$,
and suppose that $a \in K$.
Let $m\in \N$, and suppose that $\Dcoll$ satisfies the Browder condition of order $m$ at $a$, with $m$th order Browder sum $B<\infty$. Then $\delta_{a,m}$ is continuous on $R_0(K)$, with operator norm at most $B$.
\elem

Let $K$ be a compact plane set, let $a \in K$, and let $1\leq k \leq m$. Then clearly, for all $f\in R_0(K)$, we have
\[\delta_{a,k}(f)=\delta_{a,m}((z-a)^{m-k}f)\,.\]
From this is immediate that if $\delta_{a,m}$ is continuous on $R_0(K)$, then so is $\delta_{a,k}$ whenever $1\leq k \leq m$.

The following lemma is well known, but we include a sketch of the proof for the convenience of the reader.
\blem\label{EasyDerivation}
Let $K$ be a compact plane set, and let $a \in K$. Suppose that $m\in \N$ is such that $\delta_{a,m}$ is continuous on $R_0(K)$. Then there is a nondegenerate bounded point derivation of order $m$ on $R(K)$ at $a$. If $\delta_{a,m}$ is continuous on $R_0(K)$ for all $m \in \N$ then there is a nondegenerate bounded point derivation of infinite order on $R(K)$ at $a$.
\elem

\bpf
Suppose that $\delta_{a,m}$ is continuous on $R_0(K)$. Then for each $1\leq k\leq m$ we know that $\delta_{a,k}$ is also continuous on $R_0(K)$, so $\delta_{a,k}$ extends to a bounded linear functional $d^{(k)}$ on $R(K)$. Then $(d^{(k)})_{k=0}^m$ is a nondegenerate bounded point derivation of order $m$ on $R(K)$ at $a$.

The proof concerning the point derivation of infinite order is similar.
\epf

As mentioned in the introduction, Lemma \ref{Browder} is essentially immediate from \cite[Lemma 4.4]{FY}. Indeed, Lemma \ref{Browder} is an immediate corollary of Lemmas \ref{Browder2} and \ref{EasyDerivation} above.

We shall now see that the sufficient conditions in Lemma \ref{EasyDerivation} are also necessary. This result can be proved directly, but instead we will deduce it from some results for more general commutative Banach algebras.

The next result can be proved by an easy induction left to the reader.
\blem\label{EasyInduction}
Let $A$ be a commutative, unital Banach algebra, and let $\phi$ be a (nonzero) complex homomorphism on $A$. Set $M=\ker \phi$ and let $f\in M$. Let $m\in \N$, and suppose that $\der$ is a point derivation of order $m$ on $A$ at $\phi$. Then for all $k \in \N$ with $k\leq m$ we have $d^{(k)}(1)=0$, $d^{(k)}(f^k)=(d^{(1)}(f))^k$, and $d^{(k)}(M^{k+1})=0$.
\elem

The next three lemmas may be known, but we have not seen them in the literature so we prove them here.

\blem\label{codim-one}
Let $A$ be a commutative, unital Banach algebra, and let $I$ be a closed ideal in $A$.
Suppose that $\dim(I/\ol{I^2}) = 1$, and let $f \in I \sm \ol{I^2}$. Then, for all $m \in \N$, we have
\[\ol{I^m}=\ol{I^{m+1}} + \C f^m\,,\]
and so $\dim(\ol{I^m}/\ol{I^{m+1}}) \leq 1.$ In particular,
$\ol{I^m}=\ol{I^{m+1}}$ if and only if $f^m \in \ol{I^{m+1}}$.
\elem
\bpf
The last part follows immediately from the claimed equality of ideals, which we prove by induction on $m$. The result is clear when $m=1$. Now suppose that $m\in \N$ with $m>1$ and that the result is true for $m-1$. Certainly $\ol{I^{m+1}} + \C f^m$ is a closed linear subspace of $\ol{I^m}$, so it is enough to show that this subspace includes all elements of the form $gh$ where $g \in I^{m-1}$ and $h \in I$.
By the inductive hypothesis, there are scalars $\alpha$ and $\beta$, $a \in \ol{I^m}$, and $b \in \ol{I^2}$ such that $g = a+\alpha f^{m-1}$ and $h =b + \beta f.$
Then
\[gh = (a+\alpha f^{m-1})(b+\beta f) = ab + \alpha f^{m-1} b + \beta af + \alpha\beta f^m \in \ol{I^{m+1}} + \C f^m\,,\]
as required.
The induction may now proceed.
\epf

We now focus on \emph{maximal} ideals.  One easy sufficient condition for a maximal ideal $M$ to satisfy $\dim(M/\ol{M^2}) \leq 1$ is that $M$ be the closure of a principal ideal in $A$.

\blem\label{trivial-principal}
Let $A$ be a commutative, unital Banach algebra, let $M$ be a maximal ideal in $A$, let $f\in M$, and suppose that $M=\ol{fA}$. Then $\dim(M/\ol{M^2}) \leq 1$. Moreover, for all $m \in \N$, we have $\ol{M^m}=\ol{f^m A}$, and $\ol{M^m}=\ol{M^{m+1}}$ if and only if $f^m\in \ol{M^{m+1}}$.
\elem
\bpf
That $\dim(M/\ol{M^2}) \leq 1$ is immediate since $M=\ol{fA}=\C f + \ol{fM}\subset \C f + \ol{M^2}$. The rest is easy to see directly, though the last equivalence also follows from Lemma \ref{codim-one}.
\epf

We now give the connection between higher-order point derivations and powers of maximal ideals in this setting.

\blem\label{derivations-vs-ideals}
Let $A$ be a commutative, unital Banach algebra, let $\phi$ be a (nonzero) complex homomorphism
on $A$, set $M=\ker \phi$, and let $m\in\N$. Suppose that $\dim(M/\ol{M^2}) \leq 1$. Then there is a nondegenerate, bounded point derivation of order $m$ at $\phi$ on $A$ if and only if $\ol{M^m}\neq\ol{M^{m+1}}$, in which case we have
$M\supsetneq \ol{M^2} \supsetneq \cdots \supsetneq \ol{M^m} \supsetneq \ol{M^{m+1}}$.
\elem

\bpf
It is standard that if $\ol{M^2}=M$, then there are no nonzero bounded point derivations on $A$ at $\phi$. Thus we may assume that $\dim(M/\ol{M^2}) =1$.

Let $f \in M \sm \ol{M^2}\,.$ Then, by Lemma \ref{codim-one}, $\ol{M^m}\neq\ol{M^{m+1}}$ if and only if $f^m \notin \ol{M^{m+1}}$, in which case we have $f^k \in M^k \sm \ol{M^{k+1}}$ for $k=1,2,\dots,m.$

Suppose first that $\der$ is a nondegenerate bounded point derivation of order $m$ at $\phi$ on $A$. Then, by nondegeneracy, we must have $d^{(1)}(f) \neq 0$. By Lemma \ref{EasyInduction}, $d^{(m)}(f^m)=(d^{(1)}(f))^m \neq 0$, and $d^{(m)}(M^{m+1})=0$. Thus, by continuity, $f^m \notin \ol{M^{m+1}}$, as required.

Conversely, suppose that $f^m \notin \ol{M^{m+1}}$. Then an easy induction shows that
\[
A = \C 1 \oplus \C f \oplus \cdots \oplus \C f^m \oplus \ol{M^{m+1}}\,.
\]
For $0 \leq k \leq m$, we may define a bounded linear functional $d^{(k)}$ on $A$ as follows: given scalars $\alpha_j$ ($j=0,1,\dots,m$) and $g \in \ol{M^{m+1}}$, we set
\[
d^{(k)} \left(\sum_{j=0}^m \alpha_j f^j + g\right) = \alpha_k\,.
\]
It is easily verified that $\der$ is a nondegenerate bounded point derivation of order $m$ on $A$ at $\phi$.
\epf

\begin{remark}\label{trivial}
Let $K$ be a compact plane set and let $a \in K$. It is trivial to verify that, in $R(K)$, $M_a =  \ol{(z-a)R_0(K)} =\ol{(z-a)R(K)} $. It follows easily (directly, or via Lemmas \ref{codim-one} and \ref{trivial-principal}) that $\dim(M_a/\ol{M_a}^2) \leq 1$ and, for each $m \in \N$, we have $\ol{M_a^m} = \ol{(z-a)^m R(K)} = \ol{(z-a)^m R_0(K)}$. In particular, $\ol{M_a^{m+1}}=\ol{M_a^m}$ if and only if
$(z-a)^m \in \ol{(z-a)^{m+1} R_0(K)}.$
\end{remark}

\bcor\label{derivations-vs-ideals-R(X)}
Let $K$ be a compact plane set, let $a\in K$, and let $m\in\N$. Then, for the uniform algebra $R(K)$, there exists a nondegenerate bounded point derivation of order $m$ at $a$ if and only if $M_a \supsetneq \ol {M_a^2} \supsetneq \cdots \supsetneq \ol{M_a^{m+1}}$, and this holds if and only if $(z-a)^m \notin \ol{(z-a)^{m+1} R_0(K)}.$
\ecor

We are now ready to show that the sufficient conditions in
Lemma~\ref{EasyDerivation} are indeed also necessary.

\blem\label{EasyDerivation2}
Let $K$ be a compact plane set and let $a \in K$. For each $m \in \N$, there is a nondegenerate bounded point derivation of order $m$ on $R(K)$ at $a$ if and only if $\delta_{a,m}$ is continuous on $R_0(K)$. There is a nondegenerate bounded point derivation of infinite order on $R(K)$ at $a$ if and only if $\delta_{a,m}$ is continuous on $R_0(K)$ for all $m \in \N$.
\elem

\bpf
In view of Lemma \ref{EasyDerivation}, we just need to show that if $m\in \N$ and there exists a nondegenerate bounded point derivation of order $m$ on $R(K)$ at $a$, then $\delta_{a,m}$ is continuous on $R_0(K)$.
This can be proved directly for $R(K)$, or as in the proof of Lemma \ref{derivations-vs-ideals}. For the latter approach, we note that under these conditions we have
\[
R(K) = \C 1 \oplus \C (z-a) \oplus \cdots \oplus \C (z-a)^m \oplus \ol{M_a^{m+1}}\,,
\]
and then $\delta_{a,m}$ is the restriction to $R_0(K)$ of the continuous linear functional $d^{(m)}$ on $R(K)$ defined by

\[d^{(m)} \left(\sum_{j=0}^m \alpha_j (z-a)^j + g\right) = \alpha_m\,\]
where $\alpha_j$ are scalars for $j=0,1,\dots,m$ and $g \in \ol{M_a^{m+1}}$.
\epf

%
%
%
%

\section{Square roots of Swiss cheeses}\label{Square roots of Swiss cheeses}

In this section we prove Theorems~\ref{square-root-of-set} and~\ref{square-root-of-cheese}.

Recall that for a subset $E$ of $\C$, we defined the square root $\Jsqrt E$ of $E$ to be the set $\{z \in \C: z^2 \in E\,\}$. Of course the square root of a compact plane set is always a compact plane set that is symmetric about the origin.

\bpf[Proof of Theorem~\ref{square-root-of-set}]
(a) Suppose that $z^m \in \ol{J_0(R(K))}$.
Choose a sequence $\Jnseq{f_n}$ in $J_0(R(K))$ such that $f_n\ra z^m$ in $R(K)$.  Then $\Jnseq{f_n \circ z^2}$ is a sequence in $J_0(R(\Kroot))$ and $f_n\circ z^2 \ra z^{2m}$ in $R(\Kroot)$, so $z^{2m}\in \ol{J_0(R(\Kroot))}$.

(b) This is immediate from part (a) and Remark~\ref{trivial}.

(c)  Since are working with two different compact sets, we will denote the functional
$\delta_{0,m}$ on $R_0(K)$ by $\delta_{0,m}^K$ and denote the functional $\delta_{0,2m}$ on $R_0(\Kroot)$ by $\delta_{0,2m}^{\Kroot}$.

Suppose that $R(K)$ has a nondegenerate bounded point derivation of order $m$ at the origin. By Lemma \ref{EasyDerivation2}, the linear functional $\delta_{0,m}^K$ is bounded. By Lemma \ref{EasyDerivation2} again, it is sufficient to show that the linear functional $\delta_{0,2m}^{\Kroot}$ is bounded on $R_0(\Kroot)$.

Let $f \in R_0(\Kroot)$. Since $\Kroot$ is symmetric around the origin, we can set $g(z)=\bigl(f(z)+f(-z)\bigr)/2$. Then $\|g\|\leq \|f\|$ and $\delta_{0,2m}^{\Kroot}(f)=\delta_{0,2m}^{\Kroot}(g)$. Also there is $h \in R_0(K)$ with $g=h(z^2)$, and $\|h\|=\|g\| \leq \|f\|$.  (To see that $h$ really is a rational function, note that writing $f=p/q$ with $p$ and $q$ polynomials yields
\[
g(z) = \frac{1}{2}\biggl( \frac{p(z)}{q(z)} + \frac{p(-z)}{q(-z)}\biggr) =  \frac{p(z)q(-z)+q(z)p(-z)}{2q(z)q(-z)},
\]
so $g$ is a quotient of two even polynomials.)
We then have $\delta_{0,2m}^{\Kroot}(f)=\delta_{0,m}^K(h)$.

It follows that $\delta_{0,2m}^{\Kroot}$ is bounded, and indeed that $\|\delta_{0,2m}^{\Kroot}\| \leq \|\delta_{0,m}^K\|$.
\epf

To prove Theorem~\ref{square-root-of-cheese}, we need a lemma concerning square roots of open discs.

\blem\label{square-root-disc}
Let $a\in \C\setminus\{0\}$, let $0<r<|a|$, and set $\Delta=D(a,r)$, so that $r(\Delta)=r$. Set $s= s_0(\Delta)=|a|-r >0$.
Then there are two open discs $\Delta_1$ and $\Delta_2$ such that $\Jsqrt \Delta \subseteq \Delta_1 \cup \Delta_2$ and, for $i=1,2$, the following hold:
\begin{enumerate}
\item[(a)]$s_0(\Delta_i)=\sqrt{s}$;
\item[(b)]$\displaystyle r(\Delta_i)=\sqrt{|a|}-\sqrt{s} = \frac{r}{\sqrt{|a|} + \sqrt{s}} < \frac{r}{2 \sqrt s}$;
\item[(c)] for each $m \in \N$, we have
\[\frac{r(\Delta_i)}{s_0(\Delta_i)^{2m+1}} < \frac{r}{2 s^{m+1}}\,.\]
\end{enumerate}
\elem
\bpf
By rotating about the origin if necessary, we may assume that $a$ is on the positive real axis so that $|a|=a=s+r$. Note that
\begin{equation}\label{3.1.3}
\sqrt a - \sqrt s = \frac{a-s} {\sqrt a + \sqrt s} = \frac r{\sqrt a + \sqrt s}< \frac{r}{2 \sqrt s}\,.
\end{equation}
Set $a_1=\sqrt a$, $a_2=-\sqrt a$. Also set
\begin{equation}\label{3.1.4}
r_1=r_2= \sqrt a - \sqrt s = \frac{r}{\sqrt a +\sqrt s}\,,
\end{equation}
and set $\Delta_i=D(a_i,r_i)$ for $i=1,2$.
We claim that these discs satisfy the conditions of the lemma.

For condition (a) note that $s_0(\Delta_i)=\sqrt a - r_i = \sqrt s$.
Condition (b) is clear from (\ref{3.1.3}) and (\ref{3.1.4}) above.
For (c), by (a) and (b) we have
\[
\frac{r(\Delta_i)}{s_0(\Delta_i)^{2m+1}} < \frac{r}{2 \sqrt s (\sqrt s)^{2m+1}} = \frac{r}{2 s^{m+1}}\,.
\]

It remains to check the set inclusion. Let $w\in \sqrt \Delta$. Then
\[r>|w^2-a|=(w+\sqrt a)(w - \sqrt a)\,.\]
We also have $\Re(w^2)> s$, and so $|\Re(w)|>\sqrt s$.
By symmetry it is enough to consider the case $\Re(w)>0$.  Then $|w+\sqrt a|\geq \Re(w+\sqrt a) >\sqrt a + \sqrt s$, giving us
\[
|w-a_1| = |w-\sqrt a| = \frac{|w^2-a|}{|w+\sqrt a|} < \frac r{\sqrt a + \sqrt s} = r_1\,,
\]
so $w \in \Delta_1$.
\epf

\bpf[Proof of Theorem~\ref{square-root-of-cheese}]
Given a \Jgood collection $\Dcoll$ of open discs as in the statement of the theorem, for each disc $\Delta\in\Dcoll$ apply Lemma~\ref{square-root-disc} to $\Delta$ to obtain two new discs $\Delta_1$ and $\Delta_2$ as in the lemma. Then the collection $\Dtcoll$ of all new discs so obtained has all the required properties.
\epf

%
%
%
%

\section{Some preliminary results}\label{Some preliminary results}

We will need the following recent result of the second author \cite[Theorem~1.2]{Izzo} on strong regularity in $R(K)$.

\bthm\label{strongly-regular}
For each $0<r<1$, there exists a collection of open discs $\{D_k\}_{k=1}^\infty$ such that $\sum_{k=1}^\infty r(D_k)<r$ and such that, setting $K=\ol \D\sm \bigcup_{k=1}^\infty D_k$, the uniform algebra $R(K)$ is nontrivial and strongly regular.
\ethm

Recall that as mentioned in the introduction, every strongly regular uniform algebra is normal by a result of Wilken~\cite[Corollary~1]{Wilken}.
We will need the following localness property of strong regularity which is proved in \cite[Theorem~4.1]{Izzo2}.

\bthm\label{localness}
Let $A$ be a normal uniform algebra on a compact space $X$, and let $x$ be a point of $X$.  If there exists a closed neighborhood $N$ of $x$ in $X$ such that $\ol{A|N}$ (the closure in $C(N)$ of the algebra of restrictions $A|N$) is strongly regular at $x$, then $A$ is strongly regular at $x$.
\ethm

We will make repeated use of the next lemma, whose easy proof can be found in~\cite[Lemma~5.2]{Izzo}.

\blem\label{passing-to-subset}
Given compact plane sets $L\subset K$, given $m\in \N$, and given a point $a\in L$, if
$\ol{J_a(R(K))} \supset M_a(R(K))^m$, then
$\ol{J_a(R(L))} \supset M_a(R(L))^m$.
\elem

The next lemma allows us to deduce normality of $R(K)$ from normality of $R(L)$ for certain subsets $L$ of $K$.

\blem\label{easy}
Let $K$ be a compact plane set. Suppose that $R\bigl(K\sm D(0, \frac{1}{n})\bigr)$ is normal for every $n\in \N$ for which $K\sm D(0, \frac{1}{n})$ is nonempty. Then $R(K)$ is normal.
\elem

\bpf
Let $K_0$ and $K_1$ be (nonempty) disjoint closed subsets of $K$.  By reversing the roles of $K_0$ and $K_1$ if necessary, we may assume without loss of generality that the origin is \emph{not} in $K_1$.  Now choose $n\in\N$ large enough that $\ol D(0,\frac{1}{n})$ is disjoint from $K_1$.  Our hypothesis implies that there is a function $g$ in $R\bigl(K\sm D(0, \frac{1}{2n})\bigr)$ such that $g=0$ on $\bigl(K_0 \sm D(0,\frac{1}{2n})\bigr) \cup \bigl(K\cap\{\frac{1}{2n}\leq |z| \leq \frac{1}{n} \}\bigr)$, and $g=1$ on $K_1$.  The function $f$ defined on $K$ by
$$
f(z) =
\begin{cases} g(z) &\mbox{for\ } z\in K\sm D(0, \frac{1}{2n}) \\
0 & \mbox{for\ } z \in K\cap \ol D(0,\frac{1}{n})
\end{cases}
$$
is a well-defined continuous function on $K$ and is in $R(K)$ by the localization theorem for rational approximation~\cite[Theorem~26.1]{Stout}.  Since $f=0$ on $K_0$ and $f=1$ on $K_1$, this establishes the lemma.
\epf

\blem\label{iterated-cheese}
Given $v\in\No$, there exists a Swiss cheese $K$ that satisfies the Browder condition of order $2^v$ at the origin but is such that in $R(K)$ we have $\ol{J_0}\supset {M_0^{2^{v+1}}}$.
\elem

\bpf
We apply induction on $v$.  The case $v=0$ is contained in \cite[Theorem~5.1]{F0} of the first author.
Now assume that the assertion of the lemma holds for a particular value of $v$.  Then there is a \Jgood collection $\Dcoll$ of open discs that satisfies the Browder condition of order $2^v$ at the origin and is such that the set $K=\Kdef$ is a Swiss cheese as in the statement of the lemma.
Let $\Dtcoll$ be the collection of discs obtained from $\Dcoll$ in Theorem~\ref{square-root-of-cheese}, and let
$\Kt=\ol\D\sm \bigcup_{\Dt\in\Dtcoll} \Dt$.  Then $\Kt$ satisfies the Browder condition of order $2^{v+1}$  at the origin and $\ol{J_0(R(\Kt))}\supset {M_0}(R(\Kt))^{2^{v+2}}$ by Theorem~\ref{square-root-of-set}(b) and Lemma~\ref{passing-to-subset}.
\epf

A  compact plane set obtained from the closed unit disc $\ol\D$ by deleting a sequence of open discs with disjoint closures converging to the origin and whose centers lie on the positive real axis is sometimes referred to as a road runner~\cite[p.~52]{Gamelin}.
The following lemma is known, but we include a proof for completeness.

\blem\label{road-runner}
Given $m\in \N$, there exists a road runner $K$ that satisfies the Browder condition of order $m-1$ at the origin but is such that $R(K)$ does not have a nondegenerate bounded point derivation of order $m$ at the origin.
\elem

\bpf
Set $a_n= \displaystyle\frac{1}{2^n n}$ and $r_n= \displaystyle\frac{a_n^m}{2^n}$.  Set $K=\ol\D\sm \bigcup_{n=1}^\infty D(a_n,r_n)$.  Since
\[
\sum_{n=1}^\infty \frac{r_n}{(a_n-r_n)^m} \leq \sum_{m=1}^\infty \frac{r_n}{(a_n/2)^m} = \sum_{n=1}^\infty 2^{m-n} <\infty,
\]
$K$ satisfies the Browder condition of order $m-1$ at the origin.

Now consider the sequence of functions $\{f_n\}$ in $R(K)$ defined by $f_n(z)= r_n /(z-a_n)$.  Observe that $\|f_n\|_K=1$ while
$$
|\delta_{0,m}(f_n)| = r_n/ a_n^{m+1} = 1/ (2^n a_n) = n.
$$
Thus $\delta_{0,m}$ is unbounded on $R_0(K)$, and hence there is no nondegenerate bounded point derivation of order $m$ on $R(K)$ at the origin by Lemma~\ref{EasyDerivation2}.
\epf

%
%
%
%

\section{Proofs of the main theorems}\label{The Proofs of the main theorems}

\bpf[Proof of Theorem~\ref{infinite-order}]
By Theorem~\ref{strongly-regular}, for each $n\in\N$ there is a collection $\Dcoll^n$ of open discs such that $\sum_{D\in \Dcoll^n} r(D) < 1/(4n^n)$ and $R(\ol\D\sm \bigcup_{D\in\Dcoll^n} D)$ is strongly regular.
For each $n$, let $\Dtcoll^n$ be the collection of discs in $\Dcoll^n$ that intersect the annulus $K_n=\{z\in \C: \frac{1}{n} \leq |z| \leq 1\}$.

Let $K=\ol\D\sm \bigcup_n\bigcup_{\Dt^n\in\Dtcoll^n} \Dt^n$.  Then $R(K\bigcap K_n)$ is strongly regular by Lemma~\ref{passing-to-subset}, and in particular, normal.  Therefore, by Lemma~\ref{easy}, $R(K)$ is normal.  Application of Theorem~\ref{localness} now shows that $R(K)$ is strongly regular at every point other than the origin.

To complete the proof it suffices, by Lemma~\ref{Browder}, to show that the \Jgood collection $\bigcup_n \Dtcoll^n$ satisfies the Browder condition of order $m$ at the origin for every $m\in\N$.
For that observe that
\[
\sum_n \sum_{\Dt^n\in\Dtcoll^n} \frac{r(\Dt^n)}{s_0(\Dt^n)^{m+1}} \leq
\sum_n \frac{\sum_{\Dt^n\in\Dtcoll^n} r(\Dt^n)}{1/(2n)^{m+1}} \leq
 \sum_n \frac{1/(4n^n)}{1/(2n)^{m+1}}
= \sum_n \frac{2^{m-1}}{n^{n-(m+1)}}  < \infty.
\]
\epf

\bpf[Proof of Corollary~\ref{unions}]
Let $K$ be the compact plane set of Theorem~\ref{infinite-order}.  Let
$$
K_1=K\cap \left(\bigcup_{n=1}^\infty \left\{ \frac{1}{2n} \leq |z| \leq \frac{1}{2n-1}\right\} \cup \{0\}\right)
$$
and
$$
K_2=K\cap \left(\bigcup_{n=1}^\infty \left\{ \frac{1}{2n+1} \leq |z| \leq \frac{1}{2n}\right\} \cup \{0\}\right)
$$
Then $K=K_1\cup K_2$, so $R(K_1\cup K_2)$ has the properties required of it.

The strong regularity of $R(K_1\cup K_2)$ at every point other than the origin implies that each of $R(K_1)$ and $R(K_2)$ is strongly regular at every point other than the origin by Lemma~\ref{passing-to-subset}.
For strong regularity at the origin, note that for each $n\in \N$, the function that is $0$ on $K_1 \cap \bigl\{|z|< 1/(2n +\frac{1}{2})\bigr\}$ and $1$ on $K_1 \cap \bigl\{|z|> 1/(2n +\frac{1}{2})\bigr\}$ is in $R(K_1)$, and the function that is $0$ on $K_2 \cap \bigl\{|z|< 1/(2n +\frac{3}{2})\bigr\}$ and $1$ on $K_2 \cap \bigl\{|z| > 1/(2n +\frac{3}{2})\bigr\}$ is in $R(K_2)$.  Strong regularity of $R(K_1)$ and $R(K_2)$ follows.
\epf

\bpf[Proof of Theorem~\ref{main-theorem}]
Fix an integer $m\geq 2$.
Let $K_1$ be the Swiss cheese $K$ of Theorem~\ref{infinite-order}.  Let $K_2$ be the road runner $K$ of Lemma~\ref{road-runner}.  Let $K_3$ be the Swiss cheese $K$ of Lemma~\ref{iterated-cheese} with $v$ chosen large enough that $2^v +1\geq m$.  Set $K=K_1\cap K_2\cap K_3$.  Then $R(K)$ is strongly regular at every point $x\in K\sm \{0\}$ by Lemma~\ref{passing-to-subset}.  That $R(K)$ is normal follows immediately from the normality of $R(K_1)$.

Except where indicated otherwise,
it is to be understood that all ideals in this paragraph are taken in $R(K)$.
Because each of $K_1$, $K_2$, and $K_3$ satisfies the Browder condition of order $m-1$, so does $K$.  Thus $\ol{M_0^{m-1}}\supsetneq \ol{M_0^m}\supset \ol{J_0}$ by Lemmas~\ref{Browder} and~\ref{derivations-vs-ideals}.  Because $R(K_2)$ does not have a nondegenerate bounded point derivation of order $m$ at the origin, neither does $R(K)$.  Thus $\ol{M_0^m} = \ol{M_0^{m+1}}$ by Lemma~\ref{derivations-vs-ideals}, and hence $\ol{M_0^m} = \ol{M_0^q}$ for all $q\geq m$ by Lemma~\ref{straightforward}.  Finally, $\ol{J_0} \supset{M_0^{2^{v+1}}}$ because this inclusion holds for the corresponding ideals in $R(K_3)$.  Therefore, $\ol{J_0} \supset {M_0^m}$.
\epf

\bpf[Proof of Corollary~\ref{m-to-infinity}]
For each integer $m\geq2$, let $K_m$ be the Swiss cheese $K$ corresponding to $m$ in Theorem~\ref{main-theorem}.  Given a plane set $K$, a complex number $a$, and a positive real number $\rho$, denote by $a+\rho K$ the set $\{a+\rho z: z\in K\}$.  Choose a sequence of positive real numbers $\Amseq {a_m}$ monotonically strictly decreasing to 0.  Then choose a sequence of positive real numbers $\Amseq {\rho_m}$ with each $\rho_m$ small enough that the discs $a_m + \rho_m \oD$ are pairwise disjoint.  Set 
\[K
=\biggl(\,\bigcup_{m=2}^\infty (\amk) \biggr) \cup \{0\}.
\]
Application of Runge's theorem shows that 
\[
R(K)= \{ f\in C(K): f|(\amk)\in R(\amk) {\rm\ for\ all\ } m\geq 2\}.
\]
It follows easily that $R(K)$ is strongly regular at 0 and at all points of $\bigcup\limits_{m=2}^\infty (\amk)$ other than the points $a_2, a_3, \ldots$, and that for each integer $m\geq 2$ we have $\ol{J}_{a_m}\supset M_{a_m}^m$ but $\ol{J}_{a_m}\not\supset M_{a_m}^{m-1}$.

\epf

\end{document}